\documentclass[12pt]{article}
\usepackage{amsfonts,amsmath,amssymb,setspace,mathtools,bm}
\usepackage{graphicx,float}
\usepackage{amsthm}
\usepackage{fullpage}

\newtheorem{manualtheoreminner}{Theorem}
\newenvironment{manualtheorem}[1]{
  \renewcommand\themanualtheoreminner{#1}
  \manualtheoreminner
}{\endmanualtheoreminner}

\newtheorem{Theorem}{Theorem}
\newtheorem{Lemma}[Theorem]{Lemma}

\newenvironment{Proof}{\begin{trivlist} \item[] {\bf Proof.}}{\hfill $\Box$\end{trivlist}}

\renewcommand{\geq}{\geqslant}
\renewcommand{\leq}{\leqslant}

\title{The Edge-Connectivity of Vertex-Transitive Hypergraphs}

\author{Andrea C. Burgess\thanks{Department of Mathematics and Statistics, University of New Brunswick, Saint John, NB, E2L 4L5, Canada, andrea.burgess@unb.ca}
       \and Robert D. Luther\thanks{Department of Mathematics and Statistics, Memorial University of Newfoundland, St. John's, NL, A1C 5S7, Canada, rdl863@mun.ca (corresponding author)}
       \and David A. Pike\thanks{Department of Mathematics and Statistics, Memorial University of Newfoundland, St. John's, NL, A1C 5S7, Canada, dapike@mun.ca}}

\begin{document}
\date{\today}

\maketitle

\begin{abstract}
A graph or hypergraph is said to be vertex-transitive if its automorphism group acts transitively upon its vertices.
A classic theorem of Mader asserts that every connected vertex-transitive graph is maximally edge-connected.
We generalise this result to hypergraphs and show that every connected linear uniform vertex-transitive hypergraph is maximally edge-connected.
We also show that if we relax either the linear or uniform conditions in this generalisation, then we can construct examples of vertex-transitive hypergraphs which are not maximally edge-connected.
\end{abstract}


MSC Classification Codes: 05C40, 05E18, 05C65

\pagebreak
\section{Introduction}

A graph or hypergraph is {\it connected} if there is a path connecting each pair of vertices, where a {\it path} is a sequence of alternating incident vertices and edges without repetition.
A {\it cut set} of edges in a graph or hypergraph is a set of edges whose deletion renders the graph or hypergraph disconnected.
The {\it edge-connectivity} of a graph or hypergraph $H$ is the size of a minimum cut set of edges and is denoted $\kappa'(H)$.
For a graph or hypergraph $H$, $\delta(H)$ is the minimum degree among the vertices and $\Delta(H)$ is the maximum degree among the vertices, where the degree of a vertex is the number of edges incident with it.

In~\cite{W}, Whitney observes that, for a graph $G$, $\kappa'(G)$ never exceeds $\delta(G)$, a result which extends naturally to hypergraphs.
This bound is in fact tight and a graph or hypergraph $H$ which satisfies $\kappa'(H)=\delta(H)$ is said to be {\it maximally edge-connected}.
Hellwig and Volkmann list several sufficient conditions for graphs to be maximally edge-connected in their 2008 survey~\cite{HV}.

The subject of connectivity in hypergraphs has been developing recently with results like those in \cite{BS,DPP,JS}.
In~\cite{BS}, Bahmanian and \v{S}ajna study various connectivity properties in hypergraphs with an emphasis on cut sets of edges and vertices.
In~\cite{DPP}, Dewar, Pike and Proos consider both vertex and edge-connectivity in hypergraphs with additional details on the computational complexity of these problems.
In~\cite{JS}, Jami and Szigeti investigate the edge-connectivity of permutation hypergraphs.
Dankelmann and Meierling extend several well-known sufficient conditions for graphs to be maximally edge-connected to the realm of hypergraphs in~\cite{DM}.
Tong and Shan continue this work with more extensions from graphs to hypergraphs in~\cite{TS}.
Zhao and Meng present sufficient conditions for linear uniform hypergraphs to be maximally edge-connected that generalise results from graphs in~\cite{ZM}.
These three papers were primarily focused on the properties of distance and girth.

In this paper, we investigate the edge-connectivity of vertex-transitive hypergraphs which are linear and uniform.
A graph or hypergraph $H$ is said to be {\it vertex-transitive} if, for any two vertices $u$ and $v$ of $V(H)$, there exists some automorphism $\phi$ of $H$ such that $\phi(u)=v$. 
Note that any vertex-transitive graph or hypergraph must also be regular, and so $\delta(H)=\Delta(H)$.
A {\it linear} hypergraph is one in which any pair of vertices is contained in at most one edge.
A {\it uniform} hypergraph is one in which each edge has the same cardinality; moreover, if each edge has cardinality $k$, then we say that the hypergraph is $k$-uniform.

A classic result of Mader establishes the edge-connectivity of vertex-transitive graphs.
\begin{Theorem}\label{Mader}\cite{M}
Let $G$ be a vertex-transitive and connected graph. Then $G$ is maximally edge-connected.
\end{Theorem}
Our main result is a generalisation of Mader's Theorem (Theorem 1) to linear uniform hypergraphs.
In particular, we show the following:
\begin{Theorem}\label{main}
Let $H$ be a linear $k$-uniform hypergraph with $k\geq3$. If $H$ is vertex-transitive and connected, then $H$ is maximally edge-connected.
\end{Theorem}

In Section 2 we demonstrate the existence of vertex-transitive hypergraphs which fail to be maximally edge-connected when we relax either the uniformity or linearity conditions of Theorem~\ref{main}.
In Section 3 we present the proof of Theorem~\ref{main}.

\section{Non-Uniform and Non-Linear Hypergraphs}

In this section, we present two examples of vertex-transitive hypergraphs which are not maximally edge-connected. 
Both examples meet all of the criteria of the hypothesis of Theorem~\ref{main} except for linearity in the first case and uniformity in the second.

\subsection{Uniform but Non-Linear Hypergraphs}

For $k\geq 3$, let $H$ be the complete $k$-uniform hypergraph on $n$ vertices, i.e. $V(H)$ consists of $n$ vertices and $E(H)$ is equal to the set of all $k$-subsets of $V(H)$.
Then $H$ is a connected $k$-uniform hypergraph which is simple but non-linear, where a {\it simple} hypergraph is one with no repeated edges and no loops.
For any two vertices $u$ and $v$, there exists an automorphism $\phi$ such that $\phi(u)=v$, $\phi(v)=u$ and $\phi(w)=w$ for any other vertex $w$.
Therefore $H$ is also vertex-transitive.

Now let $H_1,H_2,\dots,H_k$ be distinct copies of $H$, each with its own vertex set $V(H_i)=V(H)\times \{i\}$.
Take $H^*$ to be the union of these copies along with $n$ edges of the form $E_v=\{(v,1),(v,2),\dots,(v,k)\}$ (one for each vertex $v\in V(H)$).
Then $H^*$ is a connected $k$-uniform hypergraph which is simple but non-linear.

Now we must verify that $H^*$ is vertex-transitive.
For any two vertices within the same copy of $H$, we can find an automorphism $\phi$ of $H^*$ similar to the ones described for $H$; for example, to map $(u,1)$ to $(v,1)$, use the map $\phi:H^*\rightarrow H^*$ defined by 
$$\phi(u,i)=(v,i), \phi(v,i)=(u,i) \text{ and }\phi(w,i)=(w,i) \text{ when }w\not\in\{u,v\}.$$
For any two vertices within an edge of the form $E_v$, simply take an automorphism $\psi$ of $H^*$ which exchanges the two corresponding copies of $H$ and fixes the rest; for example, to map $(v,1)$ to $(v,2)$, use the map $\psi:H^*\rightarrow H^*$ defined by 
$$\psi(u,1)=(u,2), \psi(u,2)=(u,1) \text{ and }\psi(u,i)=(u,i) \text{ when }i\not\in\{1,2\}.$$
Finally, for any two vertices in general, we may take a composition (if needed) of the two types of automorphisms we have just described.
Therefore, $H^*$ is a vertex-transitive hypergraph.
However, so long as $n \geq k+2$ and $k \geq 3$, $$\kappa'(H^*)\leq n < \binom{n-1}{k-1}+1=\Delta(H)+1=\Delta(H^*)$$ and so $H^*$ is not maximally edge-connected.

\subsection{Linear but Non-Uniform Hypergraphs}

In order to construct an example of a vertex-transitive hypergraph that is linear but non-uniform, we rely on a well known example from combinatorial designs, a finite affine plane.
A finite {\it affine plane} of order $n$ is set of $n^2+n$ lines on $n^2$ points such that each line contains $n$ points and each point lies on $n+1$ lines.
Additionally, each pair of points lie on a unique line and the lines of an affine plane can be partitioned into $n+1$ equivalence classes under the equivalence relation of parallelism; we will refer to these classes as parallel classes.
We give a direct construction of a finite affine plane of prime order as follows.

Let $k$ be an odd prime and form a $k$ by $k$ array $A$ such that the entry in row $i$ and column $j$ is $a_{i,j}=(i-1)k+j$, where $i,j\in\{1,2,\dots,k\}$.

Let the first parallel class $\Pi_0$ be the set of all rows of $A$, that is, $$\Pi_0=\big\{\{1,\dots, k\},\{k+1,\dots, 2k\},\dots,\{(k-1)k+1,\dots,k^2\}\big\}.$$
For each $i=1,\dots,k$, form the lines of parallel class $\Pi_i$ by selecting a point from row 1 and $k-1$ other points, one from each subsequent row, such that each subsequent point is located $(i-1)$ cells to the right of the last point (wrapping around if necessary).
Repeat this process for each point in row 1 to form all $k$ lines of parallel class $\Pi_i$.
Precisely, $\Pi_i$ is the collection of lines $\{B_{i,j}\}$ with $j=1,2,\dots, k$ such that each line is a set of points $B_{i,j}=\{tk+s\ |\ t\in\{0,1,\dots,k-1\} \}$ where $s$ is the unique integer between $1$ and $k$ inclusive for which $s\equiv (i-1)t+j\pmod{k}$.

Now let $H$ be the $k$-uniform hypergraph with vertex set $V(H)=\{1,2,\dots,k^2\}$ and edge set $E(H)=\bigcup_{i=1}^{k}\Pi_i$.
Note that we have intentionally left out the class $\Pi_0$.
To verify that $H$ is vertex-transitive, let $x$ and $y$ be two vertices of $H$.
Find the parallel class among $\Pi_0,\Pi_1,\dots,\Pi_k$ which contains the pair $\{x,y\}$ in a line together and write the lines of this class in order as a permutation $\sigma$.
For example, if $x$ and $y$ are both contained in the line $B_{i,j}$, then $\sigma=(L_{i,1})(L_{i,2})\cdots(L_{i,k})$, where $(L_{i,\ell})$ is a list of the points of $B_{i,\ell}$ written in a fixed order as a permutation.
Note that one of $\sigma$, $\sigma^2$,\dots,$\sigma^{k-1}$ is an automorphism in $H$ which maps $x$ to $y$ (and preserves the parallel classes).

Now take a copy of $H$ (denoted $H'$) on the vertex set $\{1',2',\dots,(k^2)'\}$ with edges corresponding to those of $H$.
Using the parallel class $\Pi_0$, form $k$ additional edges of size $2k$ as follows.
For each $i\in\{1,2,\dots,k\}$, let $e_i$ be the edge containing the $k$ vertices of the $i^{\text{th}}$ row of $A$ along with the corresponding vertices in $H'$.
In particular, for each $i\in\{1,2,\dots,k\}$,
$$e_i=\{(i-1)k+1, (i-1)k+2,\dots,(i-1)k+k, ((i-1)k+1)', ((i-1)k+2)',\dots,((i-1)k+k)'\}.$$
Then take the union of $H$, $H'$, and the $k$ edges of the form $e_i$, each of size $2k$, to form the hypergraph $H^*$.
Note that $H^*$ is a connected linear non-uniform hypergraph with edges of sizes $k$ and $2k$.
By composing the automorphisms described for $H$ with the automorphism which maps each vertex of $H$ to its copy in $H'$, we can verify that $H^*$ is also vertex-transitive.
However the edge-connectivity $\kappa'(H^*)=k$ whereas the degree $\Delta(H^*)=k+1$, and so $H^*$ is not maximally edge-connected.

\section{A Generalisation of Mader's Theorem}

Let $H$ be a hypergraph with vertex set $V(H)$. 
For $Y\subseteq V(H)$, we let $\partial(Y)$ denote the set of edges in $H$ in which each edge has at least one vertex in $Y$ and at least one vertex in $V\setminus Y$. 
A key part of the proof of our main theorem is the following lemma.

\begin{Lemma}\label{uncross}
Let $H$ be a $k$-uniform hypergraph and $X,Y\subseteq V(H)$. Then $$|\partial(X\cup Y)|+|\partial(X\cap Y)|\leq|\partial(X)|+|\partial(Y)|.$$
\end{Lemma}

\begin{Proof} In a Venn diagram of two (possibly intersecting) sets, there are four distinct regions.
For our subsets $X$ and $Y$, these are $X\setminus Y$, $Y\setminus X$, $X\cap Y$ and $(X\cup Y)^C$.
Any edges that contain vertices in more than one of these regions will contribute to the values of $|\partial(X\cup Y)|+|\partial(X\cap Y)|$ and $|\partial(X)|+|\partial(Y)|$.

When $k=2$, we have $\binom{4}{2}=6$ pairs of regions and hence, six types of relevant edges which may exist.
By checking each pair of regions, we see that $|\partial(X)|+|\partial(Y)|$ accounts for all of the edges of $|\partial(X\cup Y)|+|\partial(X\cap Y)|$ but counts any edges with vertices in both $X\setminus Y$ and $Y\setminus X$ twice, whereas $|\partial(X\cup Y)|+|\partial(X\cap Y)|$ does not count these edges at all.

When $k=3$, we have $\binom{4}{3}=4$ additional types of possible edges.
Then $|\partial(X)|+|\partial(Y)|$ accounts for all of the edges of $|\partial(X\cup Y)|+|\partial(X\cap Y)|$ but counts any edges with vertices in both $X\setminus Y$ and $Y\setminus X$ twice, whereas $|\partial(X\cup Y)|+|\partial(X\cap Y)|$ counts these edges at most once.

When $k\geq4$, there is only one additional type of possible edge, one that contains vertices from all four regions.
Such edges are contained in each of $\partial(X)$, $\partial(Y)$, $\partial(X\cup Y)$, and $\partial(X\cap Y)$, and so they are counted twice by both $|\partial(X)|+|\partial(Y)|$ and $|\partial(X\cup Y)|+|\partial(X\cap Y)|$.
\end{Proof}

We now proceed with the proof of our main result. Note that the examples detailed in Section 2 imply the necessity of the linear and uniform conditions in the statement of this result.

\begin{manualtheorem}{2}
Let $H$ be a linear $k$-uniform hypergraph with $k\geq3$. If $H$ is vertex-transitive and connected, then $H$ is maximally edge-connected.
\end{manualtheorem}

\begin{Proof} Since $\kappa'(H)\leq \Delta(H)$, it suffices to show that $\kappa'(H)\geq \Delta(H)$. Choose a proper non-empty subset $X \subset V(H)$ such that
\begin{itemize}
\item[$(i)$] $|\partial(X)|$ is minimum and 
\item[$(ii)$] $|X|$ is minimum (subject to $(i)$).
\end{itemize}
Note that by condition $(i)$, $|\partial(X)|=\kappa'(H)$, so it suffices to show that $|\partial(X)|\geq \Delta(H)$.
By definition $\partial(X)=\partial(V(H)\setminus X)$, so condition $(ii)$ implies that $|X|\leq \frac{1}{2}|V(H)|$.
In~\cite{GR} such a set $X$ is referred to as an {\it edge atom}.

Now suppose there exists $\phi \in \text{Aut}(H)$ such that $\emptyset\neq X\cap\phi(X)\neq X$.
Then by Lemma~\ref{uncross}, $$|\partial(X\cup\phi(X))|+|\partial(X\cap\phi(X))|\leq|\partial(X)|+|\partial(\phi(X))|=2|\partial(X)|.$$ 

If $|\partial(X\cup\phi(X))|< |\partial(X)|$ then the set $X\cup\phi(X)$ contradicts our choice of $X$ by condition $(i)$. 
Otherwise, $|\partial(X\cap\phi(X))|\leq|\partial(X)|$, but then $X\cap\phi(X)$ contradicts our choice of $X$ by condition $(i)$ or $(ii)$.
Therefore, for every $\phi \in \text{Aut}(H)$, either $X\cap\phi(X)=X$ or $X\cap\phi(X)=\emptyset$.
For this reason, we say that $X$ is a {\it block of imprimitivity} (for more information on this terminology, see~\cite{GR}).
This proof so far has loosely followed the proof of Mader's Theorem found in~\cite{GR}, however, to proceed from here we must make use of original techniques.

Now, for $Y\subseteq V(H)$, we let $\partial_i(Y)$ denote the set of edges in $H$ in which each edge has exactly $i$ vertices in $Y$ and $k-i$ vertices in $V\setminus Y$.
Note that $\partial(Y)=\bigcup_{i=1}^{k-1}\partial_i(Y)$.
For any $x\in X$ and $1\leq i\leq k$, let $a_i$ be the number of neighbours of $x$ in $X$ which occur in edges of $\partial_i(X)$.
Similarly, let $b_i$ be the number of neighbours of $x$ in $V\setminus X$ which occur in edges of $\partial_i(X)$.
Since $X$ is a block of imprimitivity, the values of $a_i$ and $b_i$ for $1\leq i \leq k$ do not depend on the choice of $x\in X$.

If $|X|=1$, then $|\partial(X)|=\Delta(H)$, so from now on we assume $|X|\geq 2$.
Let $x,y\in X$ and note that $$|\partial(X\setminus\{y\})|=|\partial(X)|+\frac{a_k}{k-1}-\frac{b_1}{k-1}.$$
So, if $a_k\leq b_1$ then $X\setminus\{y\}$ contradicts our choice of $X$.
Otherwise we assume $a_k>b_1$ which implies $\partial_k(X)$ is nonempty.

For the remainder of the proof, we will refer to an edge contained in the set $\partial_i(X)$ as a $\partial_i$-edge.
If $|X|=k$ then $X$ is simply a single $\partial_k$-edge.
Then by linearity, the only edges of $\partial(X)$ are $\partial_1$-edges and by vertex transitivity, $|\partial(X)|=k(\Delta(H)-1)$.
Now $|\partial(X)|=k(\Delta(H)-1)$ is strictly greater than $|\partial(\{x\})|=\Delta(H)$ as long as $k\geq 3$ and $\Delta(H)\geq 2$.
But this is easy to confirm as a connected hypergraph $H'$ with $\Delta(H')=1$ would be a single edge of $k$ vertices.
So $\{x\}$ contradicts our choice of $X$ by condition $(ii)$. Hence $|X|$ must be strictly greater than $k$.

Now since $a_k\neq 0$ and $X$ is a block of imprimitivity, every vertex of $X$ must be incident with at least one $\partial_k$-edge.
Then the collection of $\partial_k$-edges is either a collection of non-intersecting edges or a collection of edges in which each vertex of $X$ lies at the intersection of at least two of these edges.
In the first case, there must be paths in $H$ connecting the disjoint $\partial_k$-edges. But then any one of the $\partial_k$-edges would be a better choice for our set $X$ by condition $(ii)$.

Therefore, we know that each vertex of $X$ lies at the intersection of at least two $\partial_k$-edges.
For $x\in X$, let $r_x$ be the number of $\partial_k$-edges within $X$ which contain $x$.
Observe that $r_x=\frac{a_k}{k-1}$ and so $r_x$ does not depend on our choice of $x$.
So we will simply use $r$ to denote the number of $\partial_k$-edges within $X$ which contain any given vertex of $X$.
Observe that the degree of $H$, $\Delta(H)$, must be strictly greater than $r$, since otherwise every neighbour of any vertex in $X$ must also be a vertex of $X$ and therefore either $H$ is disconnected or $X=V(H)$.

In addition, we note that $\Delta(H)$ must be strictly greater than $|\partial(X)|$, since otherwise $\kappa'(H)=|\partial(X)|=\Delta(H)$.
Also $|\partial(X)|\geq \frac{|X|(\Delta(H)-r)}{k-1}$, since the edges of $\partial(X)$ can be shared by at most $k-1$ vertices of $X$.
Therefore, $$\Delta(H)>\frac{|X|(\Delta(H)-r)}{k-1};$$ rearranging for $|X|$ gives a strict upper bound $$|X|<\frac{\Delta(H)(k-1)}{\Delta(H)-r}.$$
Observe that $X$ contains the vertex $x$ and at least $r(k-1)$ other vertices.
So $$\frac{\Delta(H)(k-1)}{\Delta(H)-r}>|X|\geq 1+r(k-1).$$
This implies $\Delta(H)(k-1)>(\Delta(H)-r)+(\Delta(H)-r)r(k-1)$ and since $\Delta(H)-r>0$, we have $\Delta(H)(k-1)>(\Delta(H)-r)r(k-1)$.
Dividing both sides by $k-1\neq 0$ we have $\Delta(H)>(\Delta(H)-r)r$.

Now $\Delta(H)>(\Delta(H)-r)r$ rearranges to $r^2>\Delta(H)(r-1)$.
To make the arithmetic easier, let $d$ be the difference $\Delta(H)-r$, and note that $d$ is a positive integer.
Substitute $d+r$ for $\Delta(H)$ and continue:
$$\begin{array}{lccl}
            &           r^2 &>& (d+r)(r-1)\\
\Rightarrow &           r^2 &>& r^2+dr-d-r\\
\Rightarrow &           0   &>& dr-d-r\\
\Rightarrow &           d   &>& r(d-1).
\end{array}$$

If $d>1$ then $r<\frac{d}{d-1}$, a ratio of two consecutive positive integers, so $1\leq r<\frac{d}{d-1}\leq 2$ which implies $r=1$.
This means that each vertex of $X$ is incident with a single $\partial_k$-edge of $X$.
But we previously established that each vertex of $X$ lies at the intersection of at least two $\partial_k$-edges, a contradiction.

Finally, if $d=1$, then each vertex is incident with a single boundary edge.
Recall the lower bound $|X|\geq 1+r(k-1)$.
Replacing $r$ with $\Delta(H)-d=\Delta(H)-1$, we get $|X|\geq 1+(\Delta(H)-1)(k-1)$, which is strictly greater than $(\Delta(H)-1)(k-1)$.
So we have $\frac{|X|}{k-1}>\Delta(H)-1$.
Observe that $|\partial(X)|\geq \frac{|X|}{k-1}$ since boundary edges take up vertices of $X$ at most $k-1$ at a time.
Therefore $\frac{|X|}{k-1}>\Delta(H)-1$ implies $|\partial(X)|>\Delta(H)-1$ and so $\kappa'(H)=|\partial(X)|\geq \Delta(H)$.
\end{Proof}

\section{Acknowledgements}

Authors Burgess and Pike acknowledge NSERC Discovery Grant support and Luther acknowledges NSERC scholarship support. 


\end{document}